\documentclass{article}
\usepackage{latexsym,amssymb,graphicx}
\input amssym.def \input amssym
\newtheorem{te}{Theorem}[section]

\newtheorem{de}[te]{Definition}
\newtheorem{lm}[te]{Lemma}

\newtheorem{pp}[te]{Proposition}
\newtheorem{co}[te]{Corollary}
\newtheorem{ex}[te]{Example}
\newtheorem{qu}[te]{Problem}
\newtheorem{re}[te]{Remark}
\def\dokaz{\noindent{\bf Proof. }}
\def\kraj{\hfill $\Box$ \par \vspace*{2mm} }

\def\widemid{\hspace{1mm}\widetilde{\mid}\hspace{1mm}}

\newcommand{\zve}[1]{{{}^*\hspace{-0.5mm}#1}}

\def\gstr{\uparrow}
\def\str{\rightarrow}
\def\Str{\Rightarrow}

\def\dl{\Leftrightarrow}
\def\dstr{\hspace{-0.1cm}\downarrow}

\def\cl{{\rm cl}}

\def\cU{{\cal U}}
\def\cV{{\cal V}}
\def\cW{{\cal W}}
\def\cA{{\cal A}}

\def\cF{{\cal F}}
\def\cG{{\cal G}}
\def\cH{{\cal H}}
\def\cP{{\cal P}}
\def\cQ{{\cal Q}}

\def\bN{{\mathbb{N}}}
\def\bH{{\mathbb{H}}}

\begin{document}
\begin{center}
           {\large \bf Multiplicative finite embeddability\\
           vs divisibility of ultrafilters}\\[2mm]
{\small \bf Boris  \v Sobot}\\
{\small  Department of Mathematics and Informatics,\\
Faculty of Sciences, University of Novi Sad,\\
Trg Dositeja Obradovi\'ca 4, 21000 Novi Sad, Serbia\\
e-mail: sobot@dmi.uns.ac.rs, ORCID: 0000-0002-4848-0678}
\end{center}

\begin{abstract} \noindent
We continue the exploration of various aspects of divisibility of ultrafilters, adding one more relation to the picture: multiplicative finite embeddability. We show that it lies between divisibility relations $\mid_M$ and $\widemid$. The set of its minimal elements proves to be very rich, and the $\widemid$-hierarchy is used to get a better intuition of this richness. We find the place of the set of $\widemid$-maximal ultrafilters among some known families of ultrafilters. Finally, we introduce new notions of largeness of subsets of $\bN$, and compare it to other such notions, important for infinite combinatorics and topological dynamics.\\
%\vspace{1mm}\\


\end{abstract}

\section{Introduction}\label{intro}

As usual, by $\beta\bN$ we denote the set of ultrafilters on the set $\bN$ of natural numbers. For each $n\in\bN$, the principal ultrafilter $\{A\subseteq\bN:n\in A\}$ is identified with $n$, so $\beta\bN$ is an extension of $\bN$. The topology on $\beta\bN$ is generated by base sets $\overline{A}=\{\cF\in\beta\bN:A\in\cF\}$ for all $A\subseteq N$. All such sets are clopen. Because of nonstandard methods used in previous papers on the subject, we will avoid to denote $\overline{A}\setminus A$ by $A^*$ (so it will not get confused with a nonstandard extension $\zve A$).

Whenever $\star$ is an associative operation on $\bN$, it can be extended to an operation on $\beta\bN$ (which we also denote by $\star$):
$$\cF\star\cG=\{A\subseteq\bN:\{n\in\bN:n^{-1}A\in\cG\}\in\cF\},$$
where $n^{-1}A=\{m\in\bN:n\star m\in A\}$. Then $(\beta\bN,\star)$ becomes a compact Hausdorff right-topological semigroup. Various properties of this semigroup are described in \cite{HS}. We will be interested in extensions of operations $+$ and $\cdot$ (the usual addition and multiplication on $\bN$). Instead of $n^{-1}A$ in the case of multiplication we write $A/n=\{\frac an:a\in A\land n\mid a\}$ and in the additive case $A-n=\{a-n:a\in A\land a>n\}$. Both in $(\beta\bN,+)$ and $(\beta\bN,\cdot)$ the elements of $\bN$ are cancellative and commute with all other elements.\\

Let us fix some more notation. Caligraphic letters $\cF,\cG,\dots$ will mostly be reserved for ultrafilters, and small letters $x,y,\dots$ for elements of $\bN$. By $P$ we denote the set of prime numbers, and ultrafilters $\cP\in\overline{P}$ are called prime. Let $A$ be a subset of $\bN$. We denote $A^c=\bN\setminus A$, $A\gstr=\{n\in \bN:(\exists a\in A)a\mid n\}$, $A\dstr=\{n\in \bN:(\exists a\in A)n\mid a\}$, $\cU=\{A\in P(\bN)\setminus\{\emptyset\}:A=A\gstr\}$ and $\cV=\{A\in P(\bN)\setminus\{\bN\}:A=A\dstr\}=\{A^c:A\in\cU\}$. If $\kappa$ is a cardinal and $n\in\bN$, then $[A]^\kappa=\{X\subseteq A:|X|=\kappa\}$, $[A]^{<\kappa}=\{X\subseteq A:|X|<\kappa\}$ and $A^n=\{a^n:a\in A\}$. For $\cF\in\beta\bN$, $\cF^n$ is the ultrafilter generated by $\{A^n:A\in\cF\}$. For subsets $X$ and $Y$ of $\beta\bN$, $X\cdot Y=\{\cF\cdot\cG:\cF\in X,\cG\in Y\}$. If one of the sets $X$ or $Y$, say $X$, has only one element $\cF$, we will write $\cF\cdot Y$ instead of $\{\cF\}Y$; for example if $n\in\bN$, then $n\bN$ is the set of natural numbers divisible by $n$. We remind the reader that, for any semigroup $(\bN,*)$, $(\beta \bN,*)$ has the smallest ideal, usually denoted by $K(\beta\bN,*)$.

There are several ways to extend the divisibility relation $\mid$ from $\bN$ to $\beta\bN$. Some of them were considered in \cite{So1}:

-$\cF\mid_L\cG$ if $\cG=\cH\cdot\cF$ for some $\cH\in\beta\bN$;

-$\cF\mid_R\cG$ if $\cG=\cF\cdot\cH$ for some $\cH\in\beta\bN$;

-$\cF\mid_M\cG$ if $\cG=\cH_1\cdot\cF\cdot\cH_2$ for some $\cH_1,\cH_2\in\beta\bN$;

-$\cF\widemid\cG$ if, for every $A\in\cF$ holds $A\gstr\in\cG$.

The definition of $\widemid$ was inspired by extension of functions $f:\bN\str\bN$ to $\beta\bN$, given by: $\widetilde{f}(\cF)$ is the ultrafilter generated by $\{f[A]:A\in\cF\}$. For example, for the function $pow_n(x)=x^n$ and $\cF\in\beta\bN$, we have $\widetilde{pow_n}(\cF)=\cF^n$. The equivalent conditions for $\cF\widemid\cG$, given by
$$\cF\widemid\cG\dl\cF\cap\cU\subseteq\cG\dl\cG\cap\cV\subseteq\cF,$$
proved to be very useful. In this paper we will include another relation $\leq_{fe}$ into consideration. Strictly speaking, it was not originally defined as a divisibility relation but, as we will see in Theorem \ref{dijagram}, it is closely connected to the relations above.

Each of the mentioned relations is a quasiorder (reflexive and transitive), so by cutting with a respective equivalence relation we get orders on sets of equivalence classes. For example, we define $\cF=_\sim\cG$ if $\cF\widemid\cG$ and $\cG\widemid\cF$, so we observe $\widemid$ as an order on equivalence classes $[\cF]_\sim$ of $\beta\bN/=_\sim$. The notation for other relations is analogous.

An ultrafilter $\cF$ is divisible by $n\in\bN$ (in any of the considered relations) if and only if $n\bN\in\cF$ (\cite{So1}, Lemma 5.1). Hence we can write just $n\mid\cF$, not emphasizing the particular divisibility relation.\\

In \cite{So3} it was shown that the ultrafilters in the ``lower" part of the $\widemid$-hierarchy can be organized into $\omega$-many levels $\overline{L_n}$, where $L_n=\{p_1p_2\dots p_n:p_1,p_2,\dots,p_n\in P\}$ is the set of numbers with exactly $n$ prime (not necessarily distinct) factors. In particular $L_1=P$, so $\overline{L_1}$ is the set of all prime ultrafilters. Each proper divisor of any $\cF\in\overline{L_n}$ must belong to $\overline{L_m}$ for some $m<n$. The next result is Theorem 5.5 from \cite{So3}, with somewhat simplified notation. We remark that it will be used only once in this paper, in Theorem \ref{feminimalex}.

\begin{pp}\label{prophier}
Let $\cA_{n_1,n_2,\dots,n_m}^{\cP_1^{k_1},\cP_2^{k_2},\dots,\cP_m^{k_m}}$ be the family of all ultrafilters containing
$$F_{n_1,n_2,\dots,n_m}^{\cP_1^{k_1},\cP_2^{k_2},\dots,\cP_m^{k_m}}=\{(A_1^{k_1})^{(n_1)}(A_2^{k_2})^{(n_2)}\dots(A_m^{k_m})^{(n_m)}:$$
$$A_i\in\cP_i\land A_i\subseteq P\land (A_i=A_j\mbox{ if }\cP_i=\cP_j\mbox{ and }A_i\cap A_j=\emptyset\mbox{ otherwise})\},$$
where $\cP_1,\cP_2,\dots,\cP_m$ are prime ultrafilters, $n_1,n_2,\dots,n_m,k_1,k_2,\dots,k_m\in\bN$ and
\begin{eqnarray*}\label{eqhigher}
&& (A_1^{k_1})^{(n_1)}(A_2^{k_2})^{(n_2)}\dots(A_m^{k_m})^{(n_m)}\\
&=& \left\{\prod_{i=1}^m\prod_{j=1}^{n_i}a_{i,j}^{k_i} : a_{i,j}\in A_i\mbox{ for all }i,j\;\land\mbox{ all }a_{i,j}\mbox{ are distinct.}\right\}
\end{eqnarray*}
Then the level $\overline{L_n}$ is the disjoint union of all families $\cA_{n_1,n_2,\dots,n_m}^{\cP_1^{k_1},\cP_2^{k_2},\dots,\cP_m^{k_m}}$, for $n_1k_1+n_2k_2+\dots+n_mk_m=n$.
\end{pp}

For example, $\overline{L_2}$ consists of: (1) $\cP^2$ for $\cP\in\overline{P}$; (2) ultrafilters containing some $F_2^\cP$, and hence having only one prime $\widemid$-divisor $\cP$ and (3) ultrafilters containing $F_{1,1}^{\cP,\cQ}$ for some distinct primes $\cP$ and $\cQ$. Of course, the latter two kinds are only possible for nonprincipal $\cP$ and $\cQ$. More on this hierarchy can be found in \cite{So3}.

Ultrafilters not divisible by any $n\in\bN\setminus\{1\}$ are called $\bN$-free. Sets $A\subseteq\bN$ belonging to some $\bN$-free ultrafilter are also called $\bN$-free. Let us denote by $\cU_N$ the family of all $\bN$-free sets from $\cU$, and $\cV_N=\{A^c:A\in\cU_N\}$.\\

We recall several important notions of largeness for subsets of $\bN$. With some of them we also include equivalent conditions; proofs of these equivalences can be found in \cite{Be} and \cite{HS}. Since we want to work with two versions (multiplicative and additive) for some of the definitions, we will add the prefix A- when dealing with the $+$-version. If $\langle x_n:n\in\bN\rangle$ is a sequence of elements of $\bN$, we denote $FS(\langle x_n\rangle)=\{\sum_{n\in F}x_n:F\in[\bN]^{<\aleph_0}\}$.

A set $A\subseteq\bN$ is:

-A-IP if there is a sequence $\langle x_n:n\in \bN\rangle$ in $\bN$ such that $FS(\langle x_n\rangle)\subseteq A$, if and only if it is a member of some $+$-idempotent;%\cite{HS}, Theorem 5.12

-A-central if it is a member of some minimal $+$-idempotent (idempotent belonging to $K(\beta\bN,+)$);

-A-IP* if it intersects every A-IP set, if and only if it is a member of every $+$-idempotent;

-A-central* if it intersects every A-central set, if and only if it is a member of every minimal $+$-idempotent;

-A-thick if it contains, for every $n\in\bN$, a set $m+\{1,2,\dots,n\}=\{m+k:1\leq k\leq n\}$ for some $m\in\bN$;

-A-piecewise syndetic (short: A-pcws) if there is a finite $F\subseteq\bN$ such that the family $\{(\bigcup_{t\in F}(A-t))-n:n\in\bN\}$ has the finite intersection property, if and only if it is a member of some ultrafilter in $K(\beta\bN,+)$;%\cite{HS}, 4.40

-A-J set if for any finite set $F$ of functions from $\bN$ to $\bN$ there are $a\in\bN$ and $H\in[\bN]^{<\aleph_0}$ such that $a+\sum_{i\in H}f(i)\in A$  for all $f\in F$.\\

The corresponding multiplicative versions of these properties will be denoted with the prefix M-. For example, $A\subseteq\bN$ is M-IP if there is a sequence $\langle x_n:n\in\bN\rangle$ in $\bN$ such that $FP(\langle x_n\rangle):=\{\prod_{n\in F}x_n:F\in[\bN]^{<\aleph_0}\}\subseteq A$, if and only if it is a member of a $\cdot$-idempotent. We note also that sometimes only sets of the form $FS(\langle x_n\rangle)$ are called A-IP sets (not their supersets), but we find the definition used here more convenient.

\begin{re}\label{remdual}
In general it is common to define, for a property P, that $A\subseteq\bN$ satisfies P* if it intersects all sets with property P. Note that, if P is upwards closed (that is, $A$ has P implies that every superset of $A$ has P), then $A$ has P* if and only if $\bN\setminus A$ does not have P. It follows that, in this case, P and P* are dual: $A$ has P if and only if it intersects all sets with property P*. If we have another property Q, clearly P$\Str$Q implies Q*$\Str$P*. Hence, for upwards closed properties the reverse implication also holds.
\end{re}

The paper has six sections. Section 2 introduces multiplicative finite embeddability $\leq_{fe}$ and proves some of its basic properties. In Section 3 we explore $\leq_{fe}$-minimal ultrafilters. Using the connection between $\leq_{fe}$ and $\widemid$, we show that they are spread all over the lower half, but also present in the upper half of the $\widemid$-hierarchy. Section 4 deals with families of ultrafilters maximal for some of the relations considered. Using these families, in Section 5 we get new classes of large subsets of $\bN$ and find their place among some old notions of this kind. The closing section contains several problems we were not able to resolve.

Throughout the paper we use many results from other sources. This is natural since quite a few times we are able to draw conclusions simply  by combining these known results. Some of them we single out as separate propositions, but doing so for all of them would make the paper undesirably long.

\section{Multiplicative finite embeddability}

The concept of finite embeddability (with respect to the addition operation) was introduced (implicitly) by Ruzsa in \cite{R} and made precise by Di Nasso in \cite{D}. It was futher investigated in \cite{BD} and \cite{L2}. The paper \cite{L3} generalized the notion to $F$-finite embeddability with respect to an arbitrary set $F$ of operations. We are interested in multiplicative finite embeddability, defined with respect to the multiplication on $\bN$. Many results on $F$-finite embeddability apply to multiplicative finite embeddability, and some more results on the additive version translate directly to the multiplication version. It turns out that this relation is related to our divisibility relations. For the sake of simplicity we keep the notation $\leq_{fe}$.

\begin{de}
For $A,B\in P(\bN)$, $A$ is finitely embeddable in $B$ ($A\leq_{fe}B$) if for every finite $F\subseteq A$ there is $k\in \bN$ such that $kF\subseteq B$.

For $\cF,\cG\in\beta \bN$, $\cF$ is finitely embeddable in $\cG$ ($\cF\leq_{fe}\cG$) if for every $B\in\cG$ there is $A\in\cF$ such that $A\leq_{fe}B$.
\end{de}

%\cite{HS}, Theorem 4.47: $A\leq_{fe}B$ iff $(\zs F)(\po b\in B)bF\subseteq B$.

It is clear that the relations $\leq_{fe}$ (both on $P(\bN)$ and $\beta\bN$) are reflexive and transitive. It is a special case of \cite{L3}, Proposition 2.6 that $A\subseteq\bN$ is in the $\leq_{fe}$-greatest class of sets if and only if it is M-thick; therefore $\leq_{fe}$ for subsets of $\bN$ is not antisymmetric. Neither is $\leq_{fe}$ for ultrafilters, which will follow from Proposition \ref{femax}.

Naturally we define $\cF=_{fe}\cG$ if $\cF\leq_{fe}\cG$ and $\cG\leq_{fe}\cF$ and think of $\leq_{fe}$ as an order on equivalence classes $[\cF]_{fe}$ of $=_{fe}$. Let us also denote $\cF<_{fe}\cG$ if $\cF\leq_{fe}\cG$ but $\cF\neq_{fe}\cG$.

The additive version of the following result is included in \cite{BD}, Theorem 4 with several more equivalent conditions. We include a direct proof here.

\begin{pp}\label{equivfesets}
For $A,B\in P(\bN)$, $A\leq_{fe}B$ if and only if there is $\cH\in\beta \bN$ such that $\overline{A}\cdot\cH\subseteq\overline{B}$.
\end{pp}

\dokaz Firstly, $A\leq_{fe}B$ actually means that the family $\{B/a:a\in A\}$ has the finite intersection property. Let $\cH$ be an ultrafilter containing this family. Then, for every $\cF\in\overline{A}$, $\{n\in\bN:B/n\in\cH\}\supseteq A\in\cF$, so $B\in\cF\cdot\cH$.

In the other direction, assume that $\overline{A}\cdot\cH\subseteq\overline{B}$, but $\{B/a:a\in A\}$ does not have the finite intersection property. Then there is $a\in A$ such that $B/a\notin\cH$, so $a\cH\notin\overline{B}$, a contradiction.\kraj

%(ii) there is $\cH\in\overline{P}$ such that $A\subseteq\{n\in\bN:B/n\in\cH\}$;

%(ii)$\Str$(iii) For every $\cF\in\overline{A}$ we have $\{n\in\bN:B/n\in\cH\}\in\cF$, so $B\in\cF\cdot\cH$.

%(iii)$\Str$(i) Assume that $\overline{A}\cdot\cH\subseteq\overline{B}$, but $\{B/a:a\in A\}$ does not have the finite intersection property. Then there is $a\in A$ such that $B/a\notin\cH$, so $a\cH\notin\overline{B}$, a contradiction.\kraj

\begin{lm}\label{zaN}
$\leq_{fe}$ coincides with each of the relations $\mid_R,\mid_L,\mid_M,\widemid$ on $\bN\times\beta\bN$. Also $\cG\not\leq_{fe} m$ for all $\cG\in\beta\bN\setminus\bN$, $m\in\bN$.
\end{lm}

\dokaz As we already noted, for $\rho\in\{\mid_R,\mid_L,\mid_M,\widemid\}$, $m\in \bN$ and $\cF\in\beta \bN$, $m\rho\cF$ is equivalent to $m\bN\in\cF$. We show that the same holds for $\leq_{fe}$.

If $\cF=n\in \bN$, then $\{n\}\in\cF$ so if there are a set $A$ containing $\{m\}$ and $k\in \bN$ such that $kA\subseteq\{n\}$, it follows that $m\mid n$. On the other hand, if $m\mid n$ (say $n=mk$) then every $B\in\cF$ contains $n$, so if we take $A=\{m\}$ then $kA=\{n\}\subseteq B$.

If $\cF\in\beta \bN\setminus \bN$, assume first that $m\bN\in\cF$. Then every $B\in\cF$ contains an element $mk$ divisible by $m$, so there is $A=\{m\}$ such that $kA\subseteq B$. If we assume that $m\bN\notin\cF$, then $(m\bN)^c\in\cF$. Now there can be no set $A$ containing $\{m\}$ such that $A\leq_{fe}(m\bN)^c$: if $F=\{m\}$ then $kF\subseteq(m\bN)^c$ does not hold for any $k\in \bN$.

To prove the last statement, assume the opposite: that $\cG\leq_{fe} m$. For $B=\{m\}$ and any $A\in\cF$, if $F\subseteq A$ is of cardinality at least 2, then clearly $kF\subseteq B$ does not hold for any $k\in \bN$.\kraj

Having in mind the lemma above, in the next few proofs it suffices to consider only nonprincipal ultrafilters. The first one is a special case of \cite{L3}, Proposition 3.11, but we include a simple proof here.

\begin{lm}\label{prodfe}
If $\cF,\cG\in\beta \bN$, then $\cF\leq_{fe}\cF\cdot\cG$ and $\cG\leq_{fe}\cF\cdot\cG$.
\end{lm}

\dokaz Let $B\in\cF\cdot\cG$ be arbitrary. Then $A:=\{n\in \bN:B/n\in\cG\}\in\cF$. For every finite $F\subseteq A$ the set $\bigcap_{n\in F}B/n$ is in  $\cG$, so it is nonempty. For any $k\in\bigcap_{n\in F}B/n$ we have $kF\subseteq B$. Thus $\cF\leq_{fe}\cF\cdot\cG$.

On the other hand, for any $n\in A$ we have $B/n\in\cG$ and $n\cdot B/n\subseteq B$, so in particular for every finite $F\subseteq B/n$ we have $nF\subseteq B$. Hence $B/n\leq_{fe}B$ and $\cG\leq_{fe}\cF\cdot\cG$.\kraj

\begin{lm}\label{widefe}
$\mid_M\subset\leq_{fe}$.
\end{lm}

\dokaz From Lemma \ref{prodfe} it follows directly that $\cF\mid_R\cH$ implies $\cF\leq_{fe}\cH$ and $\cG\mid_L\cH$ implies $\cG\leq_{fe}\cH$. But $\mid_M$ is the transitive closure of $\mid_L\cup\mid_R$ and $\leq_{fe}$ is transitive, so $\mid_M\subseteq\leq_{fe}$ as well.

We will see in Proposition \ref{femax} that the $\mid_M$-maximal class is $K(\beta \bN,\cdot)$ and the $\leq_{fe}$-maximal class is $\overline{K(\beta \bN,\cdot)}$, which strictly contains $K(\beta \bN,\cdot)$. Therefore the inclusion is strict.\kraj
% Another way to conclude this (in fact, that the strict inclusion holds even at the lower levels of $L$) is shown in the following lemma.

Let us denote, for $m\in\bN$ and $A,B\in P(\bN)$, by $A\leq_{me}B$ the statement: for every $F\in[A]^m$ there is $k\in\bN$ such that $kF\subseteq B$. Then, for $\cF,\cG\in\beta\bN$, $\cF\leq_{me}\cG$ will mean: for every $B\in\cG$ there is $A\in\cF$ such that $A\leq_{me}B$.

\begin{lm}
For $\cF,\cG\in\beta\bN$, $\cF\widemid\cG$ if and only if $\cF\leq_{1e}\cG$.
\end{lm}

\dokaz If $\cF\widemid\cG$ then for any $B\in\cG$ we can take $A=B\dstr\in\cF$ and get $A\leq_{1e}B$. On the other hand, if $B\in\cG\cap\cV$, $A\leq_{1e}B$ implies $A\subseteq B$, so $A\in\cF$ implies $B\in\cF$.\kraj

\begin{lm}\label{fewide}
$\leq_{fe}\subset\widemid$.
\end{lm}

\dokaz $\leq_{fe}\subseteq\widemid$ follows directly from the definition of $\leq_{fe}$ and the lemma above. That the inclusion is strict will follow from the fact that there are (many) $\leq_{fe}$-minimal ultrafilters that are not $\widemid$-minimal (prime), see Theorems \ref{feminimal2} and \ref{feminimalex}.\kraj

The results above show that $\leq_{fe}$ is related to the divisibility relations as follows.

\begin{te}\label{dijagram}
$$\begin{array}{c}
\mid_L\\
\\
\mid_R
\end{array}\hspace{-3mm}
\begin{array}{c}
\rotatebox[origin=c]{-45}{$\subset$}\\
\rotatebox[origin=c]{45}{$\subset$}\\
\end{array}\hspace{-2mm}
\;\mid_M\;\subset\;\leq_{fe}\;\subset\widemid$$
\end{te}

The strict inclusions of the diagram above not regarding $\leq_{fe}$ were proved in \cite{So1}. Thus the following diagram is obtained, for any $\cF\in\beta\bN$:
$$\begin{array}{c}
[\cF]_L\\
\\
\;[\cF]_R\\
\end{array}\hspace{-3mm}
\begin{array}{c}
\rotatebox[origin=c]{-45}{$\subseteq$}\\
\rotatebox[origin=c]{45}{$\subseteq$}\\
\end{array}\hspace{-2mm}
\;[\cF]_M\;\subseteq\;[\cF]_{fe}\;\subseteq [\cF]_\sim$$
The inclusion between $[\cF]_M$ and $[\cF]_{fe}$ is not always equality by Proposition \ref{femax}. Neither is the inclusion between $[\cF]_{fe}$ and $[\cF]_\sim$, as we will see in Theorem \ref{kimax}.

We have a few direct corollaries of Theorem \ref{dijagram}. In a partial order with a smallest element $1$, a set $X$ is called a strong antichain if, for every two elements of $X$, $1$ is the only element below both of them.

\begin{co}
(a) For every ultrafilter $\cF\in L$, its $=_{fe}$-equivalence class is a singleton.

(b) $(\beta\bN/=_{fe},\leq_{fe})$ has strong antichains of cardinality $2^{\goth c}$.
\end{co}

\dokaz (a) By \cite{So3}, Corollary 5.10, the equivalence classes of such ultrafilters for $=_\sim$, and therefore for each of the relations $=_L$, $=_R$, $=_M$ and $=_{fe}$, are singletons.

(b) Since $(\beta\bN/=_\sim,\widemid)$ has strong antichains of cardinality $2^{\goth c}$ (for example, the antichain of all prime ultrafilters), so does $(\beta\bN/=_{fe},\leq_{fe})$.\kraj

For $p\in P$, $n\in\bN$ and $\cF\in\beta\bN$, by $p^n\parallel\cF$ we denote that $p^n\mid\cF$ but $p^{n+1}\nmid\cF$. By Lemma \ref{prodfe}, for any $\cF\in\beta\bN$ and any $p\in P$, $\cF\leq_{fe}p\cdot\cF$. It may happen that $p\cdot\cF=_{fe}\cF$, for example when $\cF$ is $\leq_{fe}$-maximal. However, if $p^n\parallel\cF$ then $\cF<_{fe}p\cdot\cF$ (because $p^{n+1}\mid\cF$, so even $\cF\neq_\sim p\cdot\cF$) and it will turn out that in this case $p\cdot\cF$ is a direct $\leq_{fe}$-successor of $\cF$.

(We know from Lemma 4.4. of \cite{So4} that $p\cdot\cF$ is indeed a direct $\widemid$-successor of $\cF$. However, this does not directly imply the same for $\leq_{fe}$: an ultrafilter $\cG$ such that $\cF<_{fe}\cG<_{fe}p\cdot\cF$ may be in the same $=_\sim$-equivalence class with $\cF$ or $p\cdot\cF$.)

\begin{te}
Let $p^n\parallel\cF$ for some $p\in P$ and some $n\in\bN\cup\{0\}$, and $\cH=p\cdot\cF$. Then there is no ultrafilter $\cG$ such that $\cF<_{fe}\cG<_{fe}\cH$.
\end{te}

\dokaz Assume the opposite. Clearly $p^n\mid\cG$. We consider two cases.

$1^\circ$ $p^{n+1}\mid\cG$, i.e.\ $p^{n+1}\bN\in\cG$. Let $B\in\cG$ be arbitrary, and let $B'=B\cap p^{n+1}\bN$. There is $A\in\cF$ such that $A\leq_{fe}B'$. Since $p^{n+1}\bN\notin\cF$, we may assume that $A\cap p^{n+1}\bN=\emptyset$. But then, for any finite $F\subseteq A$, if $kF\subseteq B'$, then $k$ must be divisible by $p$. This means that $pA\leq_{fe}B'$, and thus $pA\leq_{fe}B$. We proved that $\cH\leq_{fe}\cG$, a contradiction.

$2^\circ$ $p^{n+1}\nmid\cG$. Let $A_0\in\cF$ be arbitrary and $A:=A_0\cap p^n\bN$. Then $A\in\cF$ and $pA\in\cH$, so there is $B_0\in\cG$ such that $B_0\leq pA$. Let $B:=B_0\cap (p^n\bN\setminus p^{n+1}\bN)$. We found $B\in\cG$ such that $B\leq_{fe}A_0$, proving that $\cG\leq_{fe}\cF$, a contradiction again.\kraj

We finish this section with some results on chains in the $\leq_{fe}$-orders.

\begin{te}
In $[\bN]^{\aleph_0}$ there is a strictly $\leq_{fe}$-decreasing chain $\langle A_n:n<\omega\rangle$ without a lower bound.
\end{te}

\dokaz Let $\{F_m:m<\omega\}$ be an enumeration of $[\bN]^2$. We construct $A_n$ in such way that $F_m\not\leq_{fe}A_n$ for $m<n$. We start with $A_0:=\bN$. Assume that $A_n$ is constructed and define $A_{n+1}=\{a_m:0<m<\omega\}$ by recursion on $m$. Let $a_0$ and $a_1$ be the smallest two elements of $A_n$. Given $a_0,a_1,\dots,a_m$, for $a_{m+1}$ we choose the minimal $a\in A_n$ such that $F_0,F_1,\dots,F_n,\{a_0,a_1\}\not\leq_{fe}\{a_1,a_2,\dots,a_m,a\}$ (clearly, it suffices to take $a$ large enough). In the end, we get that $A_{n+1}\subset A_n$ (so $A_{n+1}\leq_{fe}A_n$)
 but $A_n\not\leq_{fe}A_{n+1}$. If $B\in[\bN]^{\aleph_0}$ were a lower bound for $\langle A_n:n<\omega\rangle$, for any $n$ such that $F_n\in[B]^2$ we would have $F_n\leq_{fe}A_{n+1}$, a contradiction.\kraj

What about upper and lower bounds of chains in $(\beta\bN/=_{fe},\leq_{fe})$? Using limits by ultrafilters it is not hard to construct an upper bound for any $\leq_{fe}$-chain $\langle \cF_i:i\in I\rangle$ in $\beta\bN$. Namely, if $\cV$ is any ultrafilter on $I$ containing sets $\{j\in I:j\geq i\}$ for all $i\in I$, then $\lim_{i\str\cV}\cF_i$ is the ultrafilter $\cF$ defined by: $A\in\cF$ if and only if $\{i\in I:A\in\cF_i\}\in\cV$. $\lim_{i\str\cV}\cF_i$ is an upper bound of the given chain; see Theorem 3.12 of \cite{L2} for details. In the proof of the following result bounds are also obtained as limits by ultrafilters. 

\begin{pp}\label{suplema}
(\cite{So5}, Lemma 4.1) (a) Every $\widemid$-increasing chain $\langle[\cF_i]_\sim:i\in I\rangle$ in $\beta\bN/=_\sim$ has the least upper bound $[\cG]_\sim$. If none of the ultrafilters $\cF_i$ is $\widemid$-divisible by some prime ultrafilter $\cP$, then neither is $\cG$.

(b) Every $\widemid$-decreasing chain $\langle[\cF_i]_\sim:i\in I\rangle$ in $\beta\bN/=_\sim$ has the greatest lower bound $[\cG]_\sim$.
\end{pp}

In contrast with the situation for $\widemid$-chains, for $\leq_{fe}$ we have the following result.

\begin{te}
In $\beta\bN/=_{fe}$ there is a strictly $\leq_{fe}$-increasing chain without a smallest upper bound.
\end{te}

\dokaz Let $\langle n_i:i<\omega\rangle$ be a sequence in $\bN$ such that no two differences $n_i-n_j$ for $j<i<\omega$ are the same (it is easily constructed by recursion on $i$, taking each $n_i$ large enough). Let $\cV$ be any nonprincipal ultrafilter containing $2\bN$ and $\cW$ any nonprincipal ultrafilter containing $\bN\setminus 2\bN$. Now choose any nonprincipal prime ultrafilter $\cP$ and define $\cP^{(n)}=\underbrace{\cP\cdot\cP\cdot\dots\cdot\cP}_n$. Since $\cP^{(n)}\in\overline{L_n}$, any upper bound for the sequence $\langle\cP^{(n_i)}:i<\omega\rangle$ must be in $\beta\bN\setminus L$. We consider two such upper bounds: $\cG=\lim_{i\str\cV}\cP^{(n_i)}$ and $\cH=\lim_{i\str\cW}\cP^{(n_i)}$. We will show that no upper bound for our sequence (in fact, no ultrafilter outside $L$) can be both below $\cG$ and below $\cH$.

Assume the opposite, that $\cF\in\beta\bN\setminus L$ is such that $\cF\leq_{fe}\cG$ and $\cF\leq_{fe}\cH$. If we denote $A_0=\bigcup_{i<\omega}L_{n_{2i}}$ and $A_1=\bigcup_{i<\omega}L_{n_{2i+1}}$, then $A_0\in\cG$ and $A_1\in\cH$. There are $C_0,C_1\in\cF$ such that $C_0\leq_{fe}A_0$ and $C_1\leq_{fe}A_1$. For $C=C_0\cap C_1$ we have $C\leq_{fe}A_0$ and $C\leq_{fe}A_1$. Since $\cF\notin L$, $C$ must contain elements from different levels: $c_i\in C\cap L_i$ and $c_j\in C\cap L_j$ (for some $i\neq j$). So there are $k_0,k_1\in\bN$ such that $\{k_0c_i,k_0c_j\}\subseteq A_0$ and $\{k_1c_i,k_1c_j\}\subseteq A_1$. However, if $k_0\in L_{l_0}$ and $k_1\in L_{l_1}$, then $k_0c_i\in L_{l_0+i}\cap A_0$, $k_0c_j\in L_{l_0+j}\cap A_0$, $k_1c_i\in L_{l_1+i}\cap A_1$ and $k_1c_j\in L_{l_1+j}\cap A_1$. But $(l_0+i)-(l_0+j)=(l_1+i)-(l_1+j)$, a contradiction with the choice of levels $n_i$.\kraj

\section{Minimal ultrafilters}

When speaking about minimal elements for divisibility relations, we exclude $1$, since it would trivially be the smallest for each such order. Thus, $\widemid$-minimal elements are exactly the prime ultrafilters. The set of $\mid_M$-minimal elements is the set of all ultrafilters in $\beta\bN\setminus((\beta\bN\setminus\{1\})\cdot(\beta\bN\setminus\{1\}))$. The following result says that there are many more such ultrafilters beside primes. 

\begin{pp}\label{gust}
(\cite{HS}, Theorem 6.35) $(\beta \bN\setminus \bN)\cdot(\beta \bN\setminus \bN)$ is nowhere dense in $\beta\bN\setminus\bN$, i.e.\ for every $A\in[\bN]^{\aleph_0}$ there is $X\in[A]^{\aleph_0}$ such that none of the elements of $\overline{X}$ can be represented as a product of two elements from $\beta\bN\setminus\bN$.
\end{pp}

Of course, for $\cF\in\overline{X}$ the proposition above only rules out $\mid_M$-divisors from $\beta \bN\setminus \bN$. We will see in the proof of Theorem \ref{feminimalex} how to find ultrafilters not divisible by elements of $\bN$ either.\\

It was left unresolved in \cite{L2} whether there are any minimal elements for additive finite embeddability in $\beta\bN\setminus\bN$. For multiplicative $\leq_{fe}$ we will show that there are many minimal elements in $(\beta\bN\setminus\bN,\leq_{fe})$.

\begin{lm}\label{feminimal1}
Every $\leq_{fe}$-minimal element is either a prime number $p\in P$ or an $\bN$-free ultrafilter.
\end{lm}

\dokaz Assume the opposite, that a $\leq_{fe}$-minimal ultrafilter $\cF$ is divisible by some $p\in P$, but not equal to $p$. Then $p\leq_{fe}\cF$. However, $\cF\not\neq_{fe}p$ follows from Lemma \ref{zaN}: if $\cF\in\bN$ it is a consequence of the fact that the restriction of $\leq_{fe}$ to $\bN^2$ is the ordinary divisibility, and otherwise it follows from the fact that there are no $\cF\in\beta\bN\setminus\bN$ and $n\in\bN$ such that $\cF\leq_{fe} n$.\kraj

However, not all $\bN$-free ultrafilters are $\leq_{fe}$-minimal. For example, let $\cP,\cQ\in\overline{P}\setminus P$. Then $\cP\cdot\cQ$ is $\bN$-free, but it is not minimal: by Lemma \ref{prodfe} $\cP\leq_{fe}\cP\cdot\cQ$, and $\cP\cdot\cQ\not\leq_{fe}\cP$ (because $\cP\cdot\cQ\;\widetilde{\nmid}\;\cP$).

\begin{lm}\label{nthpower}
For every ultrafilter $\cF$, $\cF^n=\cG'\cdot\cH'$ holds if and only if $\cG'=m_1\cG^n$, $\cH'=m_2\cH^n$ and $m_1m_2=m^n$ for some $m_1,m_2,m\in\bN$ and $\cG,\cH\in\beta\bN$ such that $\cF=m\cdot\cG\cdot\cH$.
\end{lm}

\dokaz First, let $\cF_1=\cG\cdot\cH$, $\cF=m\cF_1$ and $m^n=m_1m_2\in\bN^n$. For any $A\in\cF_1$ we have $B:=\{k\in\bN:A/k\in\cH\}\in\cG$. If $s=k^n\in\bN^n$, then $A^n/s=(A/k)^n$. Thus $\{s\in\bN:A^n/s\in\cH^n\}\supseteq B^n\in\cG^n$, so $A^n\in\cG^n\cdot\cH^n$. This proves that $\cF_1^n=\cG^n\cdot\cH^n$, so $\cF^n=m^n\cF_1^n=m_1\cG^n\cdot m_2\cH^n$.

In the other direction, let $\cF^n=\cG'\cdot\cH'$ for some $\cG',\cH'\in\beta\bN$. Let
$$D:=\{p_1^{s_1}p_2^{s_2}\dots p_r^{s_r}:p_1,p_2,\dots,p_r\in P\mbox{ distinct }\land 0<s_1,s_2,\dots,s_r<n\}\cup\{1\}.$$
For any $d\in\bN$, by $l(d)$ we denote the unique element od $D$ such that $d\cdot l(d)\in\bN^n$. Then $\bN^n/d=l(d)\bN^n$. Hence
$$C:=\{d\in\bN:\bN^n/d\in\cH'\}=\{d\in\bN:l(d)\bN^n\in\cH'\}$$
is nonempty only if $m_2\bN^n\in\cH'$ for some $m_2\in D$. But for distinct $d_1,d_2\in D$, $d_1\bN^n\cap d_2\bN^n=\emptyset$, so such an $m_2$ is unique. If $m_1=l(m_2)$ then $C=m_1\bN^n\in\cG'$. Hence $\cG'=m_1\cG^n$ and $\cH'=m_2\cH^n$ for some $\cG$ and $\cH$. But then $m_1m_2=m^n$ for some $m\in\bN$ and $\cF^n=m^n\cdot\cG^n\cdot\cH^n=(m\cdot\cG\cdot\cH)^n$ (by the first direction). Since $\widetilde{pow_n}$ is one-to-one, $\cF=m\cdot\cG\cdot\cH$.\kraj

\begin{te}\label{feminimal2}
For every $\leq_{fe}$-minimal ultrafilter $\cF\in\beta\bN\setminus\bN$ and every $n\in\bN$, the ultrafilter $\cF^n$ is also $\leq_{fe}$-minimal.
\end{te}

\dokaz Assume the opposite, that there is $\cG\in\beta\bN$ such that $\cG<_{fe}\cF^n$. By Lemma \ref{feminimal1} $\cF$ must be $\bN$-free, and hence so are $\cF^n$ and $\cG$. Let $B\in\cF$ be arbitrary. Then there is $A\in\cG$ such that $A\leq_{fe}B^n$. By Proposition \ref{equivfesets} there is $\cH\in\beta\bN$ such that $\overline{A}\cdot\cH\subseteq\overline{B^n}$. Then $\cG\cdot\cH=\cW^n$ for some $\cW\in\overline{B}$. By Lemma \ref{nthpower}, since $\cG$ is $\bN$-free, we must have $\cG=\cG_1^n$ for some ultrafilter $\cG_1$. But $\cG_1^n<_{fe}\cF^n$ clearly implies $\cG_1<_{fe}\cF$, a contradiction.\kraj

In particular, all prime ultrafilters must be $\leq_{fe}$-minimal since they are $\widemid$-minimal; therefore all the powers of ultrafilters from $\overline{P}\setminus P$ are also $\leq_{fe}$-minimal.

The following well-known easy fact will be useful for constructing some examples.

\begin{pp}\label{ppresek}
Every sequence $\langle X_n:n<\omega\rangle$ of infinite subsets of $\bN$ such that $X_{n+1}\subseteq X_n$ has a pseudointersection: an infinite set $Y\subseteq\bN$ such that $Y\setminus X_n$ is finite for every $n<\omega$.
\end{pp}

The next theorem shows that, in fact, there are many more $\bN$-free $\leq_{fe}$-minimal ultrafilters besides those from Theorem \ref{feminimal2}. We use notation from Proposition \ref{prophier}.

\begin{te}\label{feminimalex}
(a) Let $A_1,A_2,\dots,A_m$ be disjoint infinite subsets of $P$ and let $k_1,k_2,\dots,k_m,n_1,n_2,\dots,n_m\in\bN$ be arbitrary. Then, for some prime ultrafilters $\cP_1\in\overline{A_1},\cP_2\in\overline{A_2},\dots,\cP_m\in\overline{A_m}$, there is a $\leq_{fe}$-minimal ultrafilter $\cF\supseteq F_{n_1,n_2,\dots,n_m}^{\cP_1^{k_1},\cP_2^{k_2},\dots,\cP_m^{k_m}}$.

(b) There is a $\leq_{fe}$-minimal ultrafilter in $\beta\bN\setminus L$.
\end{te}

\dokaz (a) Let $S=(A_1^{k_1})^{(n_1)}(A_2^{k_2})^{(n_2)}\dots (A_m^{k_m})^{(n_m)}$, and let $X_n=S\cap\bigcap_{1<l\leq n}(l\bN)^c$. By Proposition \ref{ppresek} there is an infinite set $Y\subseteq S$ containing, for any $n\in\bN\setminus\{1\}$, only finitely many elements divisible by $n$. By Proposition \ref{gust} there is an infinite set $X\subseteq Y$ such that $\overline{X}\cap(\beta\bN\setminus\bN)\cdot(\beta\bN\setminus\bN)=\emptyset$. Let us show that any ultrafilter $\cF\in\overline{X}\setminus X$ is $\leq_{fe}$-minimal. Assume the opposite, that there is $\cG<_{fe}\cF$. Since $Y\in\cF$, $\cF$ is $\bN$-free, hence so is $\cG$. There is $B\in\cG$ such that $B\leq_{fe}X$, and $B$ must be infinite. By Proposition \ref{equivfesets} there is $\cH\in\beta \bN$ such that $\overline{B}\cdot\cH\subseteq\overline{X}$. Clearly $\cH\notin \bN$, since otherwise $\cH\cdot B\subseteq X\subseteq Y$. But now, for any $\cW\in\overline{B}\setminus B$, $\cW\cdot\cH\in\overline{X}\cap(\beta\bN\setminus\bN)\cdot(\beta\bN\setminus\bN)$, a contradiction. 

(b) Let $X_n=L_n\gstr\cap\bigcap_{1<l\leq n}(l\bN)^c$; by Proposition \ref{ppresek} there is a pseudointersection $Y$ of all $X_n$. Any ultrafilter $\cF\in\overline{Y}\setminus Y$ is outside $L$ and $\bN$-free. By Proposition \ref{gust} there is $X\subseteq Y$ such that all $\cF\in\overline{X}\setminus X$ are $\mid_M$-minimal. As in (a) one shows that any such $\cF$ is also $\leq_{fe}$-minimal.\kraj

Lemma 2.1 from \cite{So3} says that if, for some $A\in\cF$, $f(n)\mid n$ for all $n\in A$, then $\widetilde{f}(\cF)\widemid\cF$. Theorem \ref{feminimalex} shows that no such result holds for $\leq_{fe}$ (and consequently for $\mid_M$): if $\cF\in\overline{L_n}$ ($n>1$) is $\leq_{fe}$-minimal, any function such that $f(n)\mid n$ and $f(n)\neq n$ for all $n\in A$ (for some $A\in\cF$) would map $\cF$ to an ultrafilter $\cG\in\overline{L_m}$ for some $m<n$, but $\cG\not\leq_{fe}\cF$.

An ultrafilter $\cF$ on $\bN$ is selective (Ramsey) if, for every family $\{A_n:n<\omega\}$ of disjoint subsets of $\bN$ not in $\cF$, there is a $B\in\cF$ such that $|B\cap A_n|\leq 1$ for every $n\in\omega$. The existence of such ultrafilters can not be proved in ZFC, but if they exist, then there are also $\bN$-free selective ultrafilters: if $X$ is a pseudointersection of sets $(n\bN)^c$ for $n>1$, then any bijection $f:\bN\str X$ maps any selective ultrafilter to an $\bN$-free selective ultrafilter.

\begin{te}
For any $\bN$-free selective ultrafilter $\cF$ there is no $\cG\in L$ such that $\cG<_{fe}\cF$. In particular, if $\cF\in L$ then $\cF$ is $\leq_{fe}$-minimal.
\end{te}

\dokaz $1^\circ$ $\cF\in\overline{L_n}$ for some $n\in\omega$. Assume there is $\cG<_{fe}\cF$. $\cG$ must belong to $\overline{L_m}$ for some $m<n$. Let $k=n-m$ and let $T$ be the set of all $k$-tuples $(i_1,i_2,\dots,i_k)$ such that $1\leq i_1<i_2<\dots<i_k\leq n$. Fix a $k$-tuple $(i_1,i_2,\dots,i_k)\in T$; for every $k$-tuple $p_1\leq p_2\leq\dots\leq p_k$ of primes, let
\begin{eqnarray*}
A_{(i_1,i_2,\dots,i_k)}^{(p_1,p_2,\dots,p_k)} &=& \{q_1q_2\dots q_n:q_1,q_2,\dots,q_n\in P\land q_1\leq q_2\leq\dots\leq q_n\land\\
&& q_{i_j}=p_j\mbox{ for }j=1,2,\dots,k\}.%\mbox{ and if }q_{i_j-1}=p_j\mbox{ then }i_{j-1}=i_j-1
\end{eqnarray*}
%(The $p_j$'s are on positions $i_j$ and, if they belong to a group of successive identical factors, then they are placed at the beginning of that group.)
None of the sets $A_{(i_1,i_2,\dots,i_k)}^{(p_1,p_2,\dots,p_k)}$ is in $\cF$, because $p_1p_2\dots p_k\bN\notin\cF$. Since they are disjoint (for different $(p_1,p_2,\dots,p_k)$ and the same $(i_1,i_2,\dots,i_k)$) and $\cF$ is selective, there is $B_{(i_1,i_2,\dots,i_k)}\in\cF$ such that $|B_{(i_1,i_2,\dots,i_k)}\cap A_{(i_1,i_2,\dots,i_k)}^{(p_1,p_2,\dots,p_k)}|\leq 1$ for all $k$-tuples $p_1\leq p_2\leq\dots\leq p_k$. Hence
$$B:=\bigcap\{B_{(i_1,i_2,\dots,i_k)}:(i_1,i_2,\dots,i_k)\in T\}\in\cF$$
as well. Let $C\in\cG$ be such that $C\leq_{fe}B$; we can assume that $C\subseteq L_m$. However, if $r=|T|+1$, for any $F\in[C]^r$ there can be no $s\in\bN$ such that $sF\subseteq B$. Namely, if $s=p_1p_2\dots p_k$ where $p_1\leq p_2\leq\dots\leq p_k$, then two of the elements of $sF$ would belong to the same $A_{(i_1,i_2,\dots,i_k)}^{(p_1,p_2,\dots,p_k)}$, but $B$ intersects each such set in at most one element; a contradiction.

$2^\circ$ $\cF\notin L$. Since $\cF$ is selective, there is $B\in\cF$ such that $|B\cap L_n|\leq 1$ for all $n$. Suppose there are $m\in\bN$ and $\cG\in \overline{L_m}$ such that $\cG<_{fe}\cF$. Let $C\in\cG$ be such that $C\leq_{fe}B$; again we assume that $C\subseteq L_m$. But for any $F\in[C]^2$ and any $s\in\bN$, the elements of $sF$ must belong to the same level $L_n$, so $sF\subseteq B$ is not possible.\kraj

\section{Maximal ultrafilters}

\begin{pp}\label{femax}
(a) $K(\beta \bN,\cdot)$ is the set of $\mid_M$-maximal ultrafilters.

(b) $\overline{K(\beta \bN,\cdot)}$ is the set of $\leq_{fe}$-maximal ultrafilters.

(c) $K(\beta \bN,\cdot)$ is not topologically closed.
\end{pp}

(a) is \cite{So2}, Theorem 4.1(d), (b) is a special case of \cite{L3}, Corollary 4.10 and (c) is \cite{HS}, Corollary 8.25.

%(c) By \cite{HS}, Theorem 8.65, $\overline{K(\beta \bN,\cdot)}$ contains $2^{\goth c}$ nonminimal idempotents. By \cite{HS}, Theorem 2.9, the ideal $K(\beta \bN,\cdot)$ does not contain nonminimal idempotents.\kraj

Since $\cU$ has the finite intersection property, there is a $\widemid$-greatest class that we denote $MAX$.

\begin{pp}\label{maks}
(\cite{So4}, Lemma 4.6) For every $\cF\in\beta\bN$, $\cF\in MAX$ if and only if $n\mid\cF$ for all $n\in\bN$.
\end{pp}

The first two facts from the next lemma are a reformulation of Exercise 4.3.6 from \cite{HS}.

\begin{lm}\label{maxosob}
(a) $MAX$ is a topologically closed set.

(b) $MAX$ is a $\cdot$-ideal of $\beta \bN$.

(c) $MAX$ is closed for $+$.
\end{lm}

\dokaz (a) $MAX=\bigcap_{A\in\cU}\overline{A}$, an intersection of closed sets.

(b) Let $\cF\in MAX$ and $\cG\in\beta \bN$. Then $\cF\mid_R\cF\cdot\cG$, so by Theorem \ref{dijagram}, $\cF\widemid\cF\cdot\cG$. Hence $\cF\cdot\cG\in MAX$ as well. That $\cG\cdot\cF\in MAX$ is proved analogously.

(c) Let $\cF,\cG\in MAX$. By Proposition \ref{maks}, it suffices to show that $m\mid\cF+\cG$ for all $m\in\bN$. But $\cF=m\cF_1$ and $\cG=m\cG_1$ for some $\cF_1,\cG_1\in\beta\bN$, so $\cF+\cG=m(\cF_1+\cG_1)$ (distributivity holds in $(\beta\bN,\cdot)$ for multiplying by $m\in\bN$).\kraj

\begin{ex}
$\beta\bN\setminus MAX$ is not $\cdot$-closed; in other words we can find $\cF,\cG\notin MAX$ such that $\cF\cdot\cG\in MAX$. Let $P=Q_1\cup Q_2$ be a partition of $P$ into two infinite sets, $Q_1=\{p_n:n\in \bN\}$ and $Q_2=\{q_n:n\in \bN\}$. Let $[\cF]$ and $[\cG]$ be the smallest upper $\widemid$-bounds of $\langle\prod_{i=1}^np_i^n:n<\omega\rangle$ and $\langle\prod_{i=1}^nq_i^n:n<\omega\rangle$; here we use Proposition \ref{suplema}(a) and it is irrelevant which representatives of the $=_\sim$-equivalence classes we take. Then $\cF,\cG\notin MAX$ because $q_1\nmid\cF$ and $p_1\nmid\cG$. However, $p^i\mid\cF\cdot\cG$ for all $p\in P$ and all $i\in \bN$, so $\cF\cdot\cG\in MAX$ by Proposition \ref{maks}.
\end{ex}

\begin{lm}\label{kimax}
$\overline{K(\beta \bN,\cdot)}\subseteq MAX$.
\end{lm}

\dokaz Lemma \ref{maxosob}(b) implies that $K(\beta \bN,\cdot)\subseteq MAX$, so Lemma \ref{maxosob}(a) implies $\overline{K(\beta \bN,\cdot)}\subseteq MAX$.\kraj

The strict inclusion will be proved in Corollary \ref{cormax}.

An ultrafilter $\cF\in\beta\bN$ is called weakly summable if every $A\in\cF$ is an A-IP set. In \cite{HS}, Theorem 12.17, it was shown that, if we denote $\Gamma=\cl\{\cF\in\beta \bN:\cF+\cF=\cF\}$, then $\Gamma$ is exactly the set of weakly summable ultrafilters. Thus every $+$-idempotent ultrafilter is weakly summable; the reverse is false by \cite{HS}, Theorem 8.24 or Theorem 12.22. The following proposition is (a) \cite{HS}, Theorem 5.20 or \cite{Be}, Lemma 4.1, and (b) \cite{HS}, Theorem 12.23.

\begin{pp}\label{gammaideal}
(a) $\Gamma$ is a left $\cdot$-ideal of $\beta\bN$;

(b) $\Gamma$ is not $+$-closed.
\end{pp}

\begin{te}\label{gamma}
$\Gamma\subset MAX$.
\end{te}

\dokaz By \cite{HS}, Theorem 5.19.1 every $+$-idempotent is divisible by all $n\in \bN$, so Proposition \ref{maks} implies that the set of all $+$-idempotents is included in $MAX$. By Lemma \ref{maxosob}(a) MAX also contains its closure $\Gamma$.

However, $\Gamma\neq MAX$. One way to see this is the fact that $MAX$ is $+$-closed and $\Gamma$ is not. The other way is described in Example \ref{exgamma} below.\kraj

\begin{ex}\label{exgamma}
Let $\langle a_n:n\in \bN\rangle$ be a sequence in $\bN$ defined recursively: $a_1=1$ and, for all $n>1$, $n\mid a_n$ and $a_n>\sum_{i=1}^{n-1}a_i$. Then the set $A=\{a_n:n\in \bN\}$ intersects $m\bN$ for every $m\in \bN$. Thus $\cU\cup\{A\}$ has the finite intersection property, and there is $\cF\in MAX$ such that $A\in\cF$. But $A$ is not an A-IP set, so there is no $+$-idempotent containing $A$, and $\cF\notin\Gamma$.
\end{ex}

The set $\bH=\bigcap_{n\in\bN}\overline{2^n\bN}$ plays an important role in the investigation of the algebraic structure of $\beta\bN$, see Chapter 6 of \cite{HS}. It is obvious that $MAX\subset\bH$, and by \cite{HS}, Corollary 7.26, $MAX$ is algebraically and topologically isomorphic to $\bH$. The results above give us the following image, showing the place of $MAX$ among these important families of ultrafilters.

\begin{center}
\includegraphics[scale=0.07]{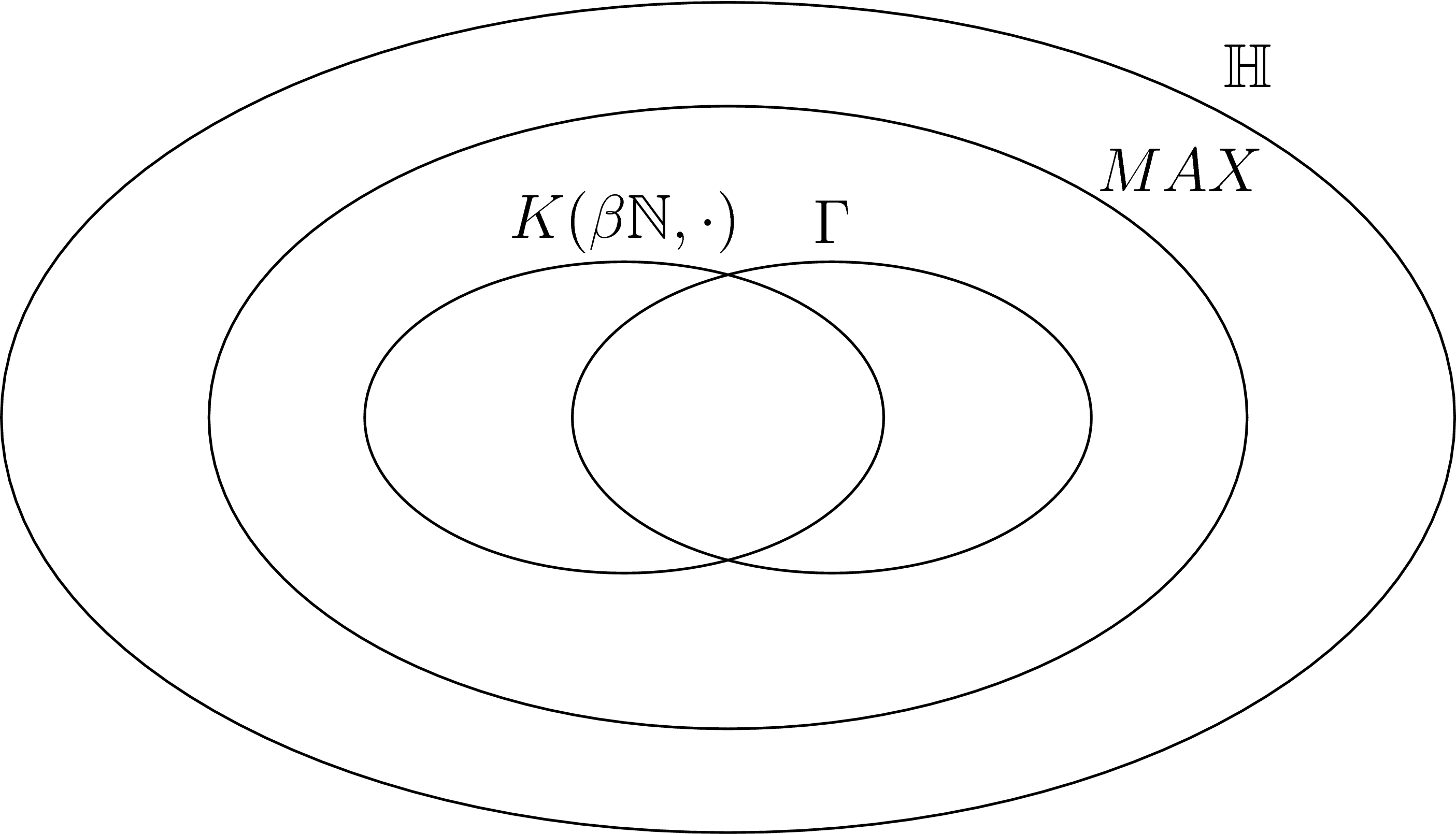}
\end{center}

Recall that $\cU_N$ is the family of all $\bN$-free sets from $\cU$, and $\cV_N=\{A^c:A\in\cU_N\}$. In \cite{So5} it was shown that $\cU_N\cup\{(n\bN)^c:n\in\bN\setminus\{1\}\}$ has the finite intersection property, so there is a $\widemid$-greatest $\bN$-free class that we denote $NMAX$. $A\subseteq\bN$ is a strong antichain if every two elements $m,n\in A$ are mutually prime.

\begin{pp}\label{nmax}
(\cite{So5}, Lemma 5.2) $S\subseteq\bN$ is an $\bN$-free set if and only if it contains an infinite strong antichain.
\end{pp}

\begin{lm}\label{nmaxosob}
(a) $NMAX$ is a topologically closed set.

(b) $NMAX$ is a $\cdot$-ideal of $\beta\bN_F:=\{\cG\in\beta\bN:\cG\mbox{ is }\bN$-free$\}$.
\end{lm}

\dokaz (a) $NMAX=\bigcap_{A\in\cU_N}\overline{A}\cap\bigcap_{n\in\bN\setminus\{1\}}\overline{(n\bN)^c}$, an intersection of closed sets.

(b) Let $\cF\in NMAX$ and $\cG\in\beta\bN_F$. Assume $p\mid\cF\cdot\cG$ for some $p\in P$. Then $p\mid\cF$ or $p\mid\cG$, a contradiction. So $\cF\cdot\cG$ is $\bN$-free. But it is $\mid_M$-divisible by $\cF$ and hence $\widemid$-divisible by all $\bN$-free ultrafilters, so $\cF\cdot\cG\in NMAX$.\kraj

$NMAX$ is clearly not an ideal in $\beta\bN$, since for $\cF\in NMAX$ and $n\in\bN\setminus\{1\}$ we have $n\cdot\cF\notin NMAX$. Also, it is not $+$-closed because a sum of two ultrafilters $\cF,\cG$ not divisible by 2 must be divisible by 2: if we denote by $2\bN+1$ the set of odd numbers, then $\{n\in\bN:(2\bN+1)-n\in\cG\}=2\bN\in\cF$, so $2\bN\in\cF+\cG$.

\section{Large subsets of $\bN$}

To definitions of various largeness properties from Section \ref{intro} we now add a few more, derived from our divisibility relation $\widemid$. For each of them we also give some immediate equivalents. A set $A\subseteq\bN$ is:

-MAX-set if it a member of some ultrafilter in $MAX$, if and only if it is not contained in a set in $\cV$, if and only if $A\dstr=\bN$;

-NMAX-set if it is a member of an ultrafilter divisible by all $\bN$-free ultrafilters, if and only if it is not contained in a set in $\cV_N$;

-MAX*-set if it intersects every MAX-set, if and only if it is a member of every ultrafilter in $MAX$, if and only if it contains some set in $\cU$;

-NMAX*-set if it intersects every NMAX-set, if and only if it is a member of every ultrafilter divisible by all $\bN$-free ultrafilters, if and only if it contains some set in $\cU_N$.

We will obtain the following situation:
\begin{center}
\includegraphics[scale=0.7]{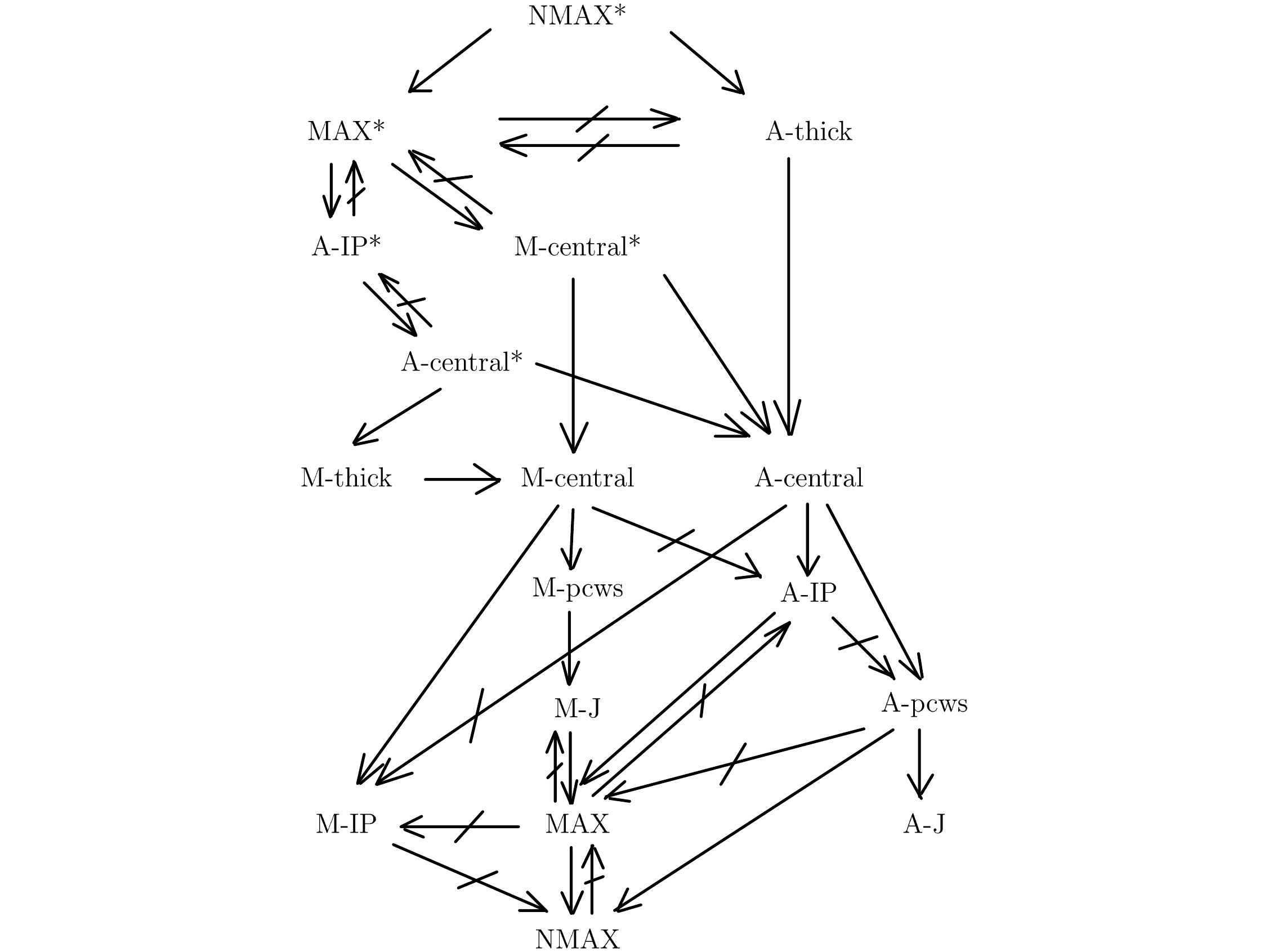}
\end{center}

Some of the implications are trivial: MAX$\Str$NMAX, NMAX$\not\Str$MAX, NMAX* $\Str$MAX*, A-IP*$\Str$A-central*$\Str$A-central$\Str$A-IP and M-central*$\Str$M-central$\Str$M-IP. Proofs of most of the other known implications can be found in \cite{Be}: A-central*$\Str$M-thick is 5.11, A-thick$\Str$A-central and M-thick$\Str$M-central are 5.10, A-central$\Str$A-pcws is 5.5, M-central$\Str$M-pcws is 5.7, M-central$\not\Str$A-IP is 5.16 and A-central$\not\Str$M-IP is Proposition 5.14. A-central*$\not\Str$A-IP*, M-central*$\Str$A-central and A-pcws$\Str$A-J (as well as M-pcws$\Str$M-J) are, respectively, Theorem 16.8, Theorem 16.26.1 and Theorem 14.8.3 from \cite{HS}. For sake of completeness, we include the proof of A-IP$\not\Str$A-pcws in the theorem below.

Now let us consider the remaining implications, involving sets derived from maximal divisibility classes.

\begin{te}\label{veliki}
For subsets of $\bN$ the following implications hold:

(a) {\rm A-IP}$\Str${\rm MAX} but {\rm MAX}$\not\Str${\rm A-IP}.

(b) {\rm M-IP}$\not\Str${\rm NMAX}.

(c) {\rm M-J}$\Str${\rm MAX} but {\rm MAX}$\not\Str${\rm M-J}.

(d) {\rm A-pcws}$\not\Str${\rm MAX}.

(e) {\rm A-pcws}$\Str${\rm NMAX}.

(f) {\rm A-IP}$\not\Str${\rm A-pcws}.
\end{te}

\dokaz (a) Every A-IP set is contained in a $+$-idempotent, which is by Lemma \ref{gamma} an element of $MAX$. On the other hand, in Example \ref{exgamma} we actually constructed a set belonging to an ultrafilter in $MAX$ which is not an A-IP set.

%(b) Similar to MAX$\not\Str$A-IP. We construct elements of a set $A=\{a_n:n\in\bN\}$ recursively: $a_1=1$ and, for all $n>1$, $n\mid a_n$ and $a_n>\prod_{i=1}^{n-1}a_i$. Then $A$ is not an M-IP set but $\cU\cup\{A\}$ has the finite intersection property, so $A$ is an element of an ultrafilter in $MAX$.

(b) Let $P=\{p_n:n<\omega\}\cup Q$ be a partition of the set $P$ of primes into two infinite subsets. If $A=FP(\langle p_n\rangle)$, then $A$ is an M-IP set. However, by Proposition \ref{nmax} $Q\gstr$ is a set in $\cU_N$ and $A$ is contained in its complement, which is in $\cV_N$.

(c) Assume the opposite, that there is a M-J-set $A$ such that $A\subseteq B$ for some $B\in\cV$. Any number divisible by any element of the complement $B^c$ also belongs to $B^c$. Hence, if we define a function $f:\bN\str B^c$ arbitrarily, $a\cdot\prod_{i\in H}f(i)\notin B$ for all $a\in\bN$, $H\in[\bN]^{<\aleph_0}$, a contradiction.

To prove that the reverse implication does not hold, let $P=\{p_n:n\in\bN\}$ be the increasing enumeration of the primes and
$$A=\{p_{i_1}^kp_{i_2}^k\dots p_{i_m}^k:m,k\in\bN\land i_1<i_2<\dots<i_m\}.$$
Then $A\dstr=\bN$, so $A$ is a MAX-set. To prove that it is not M-J, let functions $f,g:\bN\str\bN$ be defined by $f(n)=p_{2n}^2p_{2n+1}$ and $g(n)=p_{2n}p_{2n+1}^2$. Now, if $H\in[\bN]^{<\aleph_0}$, $a\in\bN$, $x_f=a\prod_{i\in H}f(i)$ and $x_g=a\prod_{i\in H}g(i)$, let $j=\min H$, and $p_{2j}^k\parallel a$ and $p_{2j+1}^l\parallel a$ (i.e.\ $k,l$ are the largest numbers for which these divisibilities hold). Then $x_f=p_{2j}^{k+2}p_{2j+1}^{l+1}q_f$ and $x_g=p_{2j}^{k+1}p_{2j+1}^{l+2}q_g$ for some $q_f,q_g$ not divisible by $p_{2j}$ or $p_{2j+1}$. But it is impossible that both $k+2=l+1$ and $k+1=l+2$ hold, so $x_f$ and $x_g$ can not both belong to $A$.

(d) The set $2\bN+1$ of odd numbers is A-piecewise syndetic, but is clearly not a MAX-set.

(e) Let $A$ be a set not in NMAX; we may assume $A\in\cV_N$. Suppose that $A$ is A-pcws. Then there is finite $F=\{t_1,t_2,\dots,t_f\}\subset\bN$ such that, for $B=\bigcup_{t\in F}(A-t)$, the set $\{B-k:k\in\bN\}$ has the finite intersection property. Since $A^c$ is $\bN$-free, by Proposition \ref{nmax} there are mutually prime $n_1,n_2,\dots,n_f\in A^c$. Let $n=n_1n_2\dots n_f$; then there is $m\in\bigcap_{1\leq i\leq n}(B-i)$. This means that $\{m+1,m+2,\dots,m+n\}\subseteq B$. By the Chinese remainder theorem, the system
$$x\equiv -t_j(\mbox{mod }n_j),\;\;j=1,2,\dots,f$$
has a solution $x\in\{m+1,m+2,\dots,m+n\}$. Then $n_j\mid x+t_j$ and $n_j\in A^c$ imply $x+t_j\in A^c$, since $A^c\in\cU_N$. However, since $x\in B$, $x+t_j\in A$ for some $1\leq j\leq f$; a contradiction.

(f) Let $\langle a_n:n\in\bN\rangle$ be a sequence in $\bN$ such that $a_n>n+\sum_{i=1}^{n-1}a_i$ for $n>1$. Then $A:=FS(\langle a_n\rangle)$ is an A-IP set. Assume that, for some $n\in\bN$ and $B=\bigcup_{t\leq n}(A-t)$, the family $\{B-k:k\in\bN\}$ has the finite intersection property. Then there is $x\in\bigcap_{k\leq 2^{n-1}}(B-k)$, which means that $B$ contains $2^{n-1}$ successive numbers. But each of them is a nonempty sum of some of $a_n$ minus some $t$. Hence at least $n$ elements of $\{a_n:n\in\bN\}$ must be "used" to get these elements, say $a_{i_1}<a_{i_2}<\dots<a_{i_n}<\dots$ If $y$ is the smallest among them which contains $a_{i_m}$ for some $m\geq n$, say $y=a_{i_m}+a_{j_1}+\dots+a_{j_l}-t$ for some $t\leq n$, then $y-1$ can not belong to $B$, since any element of that set not using $a_{i_k}$ for $k\geq m$ is no larger than $a_1+a_2+\dots+a_{i_m-1}<a_{i_m}-i_m\leq y-1$.\kraj

\begin{te}\label{veliki2}
For subsets of $\bN$ the following implications hold:

(a) {\rm MAX*}$\Str${\rm M-central*} but {\rm M-central*}$\not\Str${\rm MAX*}.

(b) {\rm MAX*}$\Str${\rm A-IP*} but {\rm A-IP*}$\not\Str${\rm MAX*}.

(c) {\rm A-thick}$\not\Str${\rm MAX*} and {\rm MAX*}$\not\Str${\rm A-thick}.

(d) {\rm NMAX*}$\Str${\rm A-thick}.
\end{te}

\dokaz (a) We use Remark \ref{remdual}. M-central$\Str$MAX directly implies that MAX*$\Str$M-central*. If M-central*$\Str$MAX* were true, by duality we would have MAX$\Str$M-central, which is not true.

(b) By Remark \ref{remdual}, MAX*$\Str$A-IP* follows from Theorem \ref{veliki}(a). We can use duality to prove {\rm A-IP*}$\not\Str${\rm MAX*} as well, but here is a direct proof. Since $\Gamma$ is not $+$-closed (Proposition \ref{gammaideal}), there are $\cF,\cG\in\Gamma$ such that $\cF+\cG\notin\Gamma$. By Lemmas \ref{gamma} and \ref{maxosob}(c), $\cF+\cG\in MAX$. $\Gamma$ is topologically closed, so there is $A\in\cF+\cG$ such that $\overline{A}\cap\Gamma=\emptyset$. This means that $A^c$ is an A-IP* set, but not a MAX*-set.

(c) We construct finite sets $F_n$ and elements $a_n\in\bN$ by recursion on $n$. Let $F_0=\emptyset$ and $a_0=1$. For $n>1$, let $F_n$ be any set of $n$ successive natural numbers larger than $a_{n-1}$, and let $a_n$ be larger than all elements of $F_n$ and divisible by $n$. Then the set $\bigcup_{n<\omega}F_n$ is A-thick but not in MAX* (it is disjoint from $\{a_n:n<\omega\}$, so it does not contain a set in $\cU$).

To show that MAX*$\not\Str$A-thick, note just that $2\bN$ is a MAX*-set but does not even contain an interval of length 2.

(d) If $X$ is a NMAX*-set, we can assume that $X\in\cU_N$. By Proposition \ref{nmax} there is an infinite strong antichain $A=\{a_m:m\in\bN\}\subseteq X$. Let $n\in\bN$ be given; we are looking for an interval of length $n$ in $X$. Consider the system of congruences
$$x\equiv -m(\mbox{mod }a_m),\;\;m=1,2,\dots,n.$$
Since $a_m$ are pairwise mutually prime, this system has a solution $x\in\bN$ by the Chinese remainder theorem. Then, for each $m\in\{1,2,\dots,n\}$, $x+m$ is divisible by $a_m$, so $a_m\in X$ implies $x+m\in X$.\kraj

There are several direct corollaries of the results above.

\begin{co}\label{corcor}
(a) {\rm A-pcws}$\not\Str${\rm A-IP}.

(b) {\rm MAX}$\not\Str${\rm A-pcws}.

(c) {\rm MAX}$\not\Str${\rm M-pcws}.

(d) {\rm M-central*}$\not\Str${\rm A-thick}.

(e) {\rm MAX}$\not\Str${\rm M-IP}.

(f) {\rm M-IP}$\not\Str${\rm M-J}.
\end{co}

\begin{co}\label{cormax}
(a) $\overline{K(\beta \bN,\cdot)}\subset MAX$.

(b) $K(\beta\bN,+)\not\subseteq MAX$.
\end{co}

\dokaz (a) By Lemma \ref{kimax} we only need to show that the inclusion is strict. By Corollary \ref{corcor}(c) there is a set $A$ belonging to an ultrafilter $\cF\in MAX$ which is not M-piecewise syndetic. This means that $A$ is not an element of any ultrafilter in $K(\beta\bN,\cdot)$. Since $\overline{A}$ is an open set, it is disjoint from $\overline{K(\beta\bN,\cdot)}$, so $\cF\notin\overline{K(\beta\bN,\cdot)}$.

(b) By Theorem \ref{veliki}(d) there is $A\subseteq\bN$ which is A-piecewise syndetic set but not a MAX-set. Then there is an ultrafilter $\cF\in\overline{A}\cap K(\beta\bN,+)$, but there is no ultrafilter from $MAX$ in $\overline{A}$.\kraj

\section{Open problems and remarks}

The $\leq_{fe}$ relation for sets has a nice nonstandard characterization: by (an analogon of) \cite{D}, Proposition 4.3, $A\leq_{fe}B$ holds if and only if $x\cdot A\subseteq\zve B$ for some $x\in\zve\bN$. It would be very useful if such a result could be found for the ultrafilter relation, but we were not able to find it so far.

\begin{qu}
Find a nice nonstandard characterization of $\leq_{fe}$ for ultrafilters.
\end{qu}

The only examples of ultrafilters such that $\cF\leq_{fe}\cG$ but $\cF\nmid_M\cG$ we found so far are above $L$. So it is natural to ask the following.

\begin{qu}
Does $\cF\leq_{fe}\cG$ imply $\cF\mid_M\cG$ for $\cF,\cG\in L$?
\end{qu}

\begin{qu}\label{probmin}
(a) Find a nice characterization of $\leq_{fe}$-minimal ultrafilters.

(b) Is every $\mid_M$-minimal ultrafilter also $\leq_{fe}$-minimal?
\end{qu}

If the answer to (b) is yes, then a solution to (a) would also provide a nice characterization of $\mid_M$-minimal ultrafilters.

In Section 2 we considered upper bounds of chains in $(\beta\bN/=_{fe},\leq_{fe})$, but did not resolve the analogous questions for lower bounds.

\begin{qu}
Does every $\leq_{fe}$-decreasing chain in $\beta\bN/=_{fe}$ has a lower bound? If yes, does it always have a greatest lower bound?
\end{qu}

Concerning Corollary \ref{cormax} and Theorem \ref{gamma} we have the following question.

\begin{qu}
(a) Is MAX the $\cdot$-ideal of $\beta \bN$ generated by $\Gamma$?

(b) Is MAX the $+$-closure of $\Gamma$ or the $+$-closure of $\overline{K(\beta \bN,\cdot)}$?
\end{qu}

Finally, the last problem deals with a missing implication on the diagram from Section 5.

\begin{qu}
Does A-J$\Str$NMAX?
\end{qu}

The author gratefully acknowledges financial support of the Science Fund of the Republic of Serbia (call PROMIS, project CLOUDS, grant no.\ 6062228) and Ministry of Education, Science and Technological Development of the Republic of Serbia (grant no.\ 451-03-9/2021-14/200125).

\end{document}